\documentclass[10pt]{amsart}
\usepackage{amsmath,amssymb,latexsym,hyperref,url}

\newcommand{\seqnum}[1]{\href{http://www.research.att.com/~njas/sequences/#1}{\underline{#1}}}

\begin{document}

\title{On the number of walks on a regular Cayley tree}
\author{Eric Rowland}
%\email{erowland@math.rutgers.edu}
\author{Doron Zeilberger}
%\email{zeilberg@math.rutgers.edu}
%\address{
%Department of Mathematics, Rutgers University (New Brunswick),
%Hill Center-Busch Campus, 110 Frelinghuysen Rd., Piscataway,
%NJ 08854-8019, USA.
%}
\date{March 12, 2009}

\begin{abstract}
We provide a new 
derivation of the 
well-known generating function counting the 
number of walks on a regular tree that start and end at the same vertex,
and more generally, a generating function for the number of walks 
that end at a vertex a distance $i$ from the start vertex.
These formulas seem to be very old, and  go back, in an equivalent form,
at least to Harry Kesten's work on symmetric random walks on groups
from 1959, and in the present form to Brendan McKay (1983).
\end{abstract}

\maketitle
%\markboth{}{}

Consider the (infinite) regular Cayley tree where each vertex has degree $m$.
How many walks of length $n$ are there that start and end at the
same vertex? Let's call this number $A_m(0,n)$. More generally,
how many such walks of length $n$ are there that end up at a given vertex that has distance
$i$ from the starting vertex? Let's call that number $A_m(i,n)$.
The sequences $\{A_m(0,2n)\}_{n=0}^{\infty}$ for $m=2,3,4$ are Sloane's
\seqnum{A000984} (the central binomial coefficients),
\seqnum{A089022}, and \seqnum{A035610} respectively \cite{eis}. For $5 \leq m \leq 8$ these
are sequences $\text{A}12096r, 0 \leq r \leq 3$. The sequences for $i>0$ do not
seem to be present at the time of this writing.

The generating function of $\{A_m(0,n)\}$,
$$
f_m(t):=\sum_{n=0}^{\infty} A_m(0,n)t^{n},
$$
supplied to Sloane \cite[\seqnum{A035610}]{eis} by Paul Boddington, is
$$
f_m(t)=\frac{2(m-1)}{m-2+m\sqrt{1-4(m-1)t^2}}.
$$
We couldn't find any proof of this in the literature.  Stevanovi\'c et al. \cite{stevanovic} reference Boddington's formula without proof.
After the first version of this
article was written, Brendan McKay pointed out to us
that the formula for $f_m(t)$, as well as
the formula for $f_m^{(i)}(t)$ given below,
are contained in \cite{mckay}. Shortly after, we
also received an email message from Franz Lehner, who
kindly told us about Harry Kesten's work \cite{kesten}.
He also mentioned Pierre Cartier's work \cite{cartier}.

Nevertheless, we still believe that the present note
is worth publishing, because of the elegant method of proof
that should be applicable in many other problems.
Of course, none of this is ``deep'' (by today's state of the art),
so our paper should be labeled ``for entertainment only''.

More generally, we prove that
$$
f_m^{(i)}(t):=\sum_{n=0}^{\infty} A_m(i,n)t^n,
$$
the generating function for those walks that end up at a vertex that has distance $i$ from the starting vertex, for any $i \geq 0$,
equals
$$
f_m^{(i)}(t)=
\frac{2(m-1)}{m-2+m\sqrt{1-4(m-1)t^2}} \cdot
\left( \frac{ 1 - \sqrt{1-4(m-1)t^2} }{2(m-1)t} \right )^i.
$$
Let us remark that $A_{2m}(i,n)$ is
also the coefficient of any reduced word of length $i$ in the 
expansion of $\left( \sum_{i=1}^{m} x_i+x_i^{-1} \right)^n$ in non-commuting symbols $x_1, \dots, x_m$.

\section*{The recurrence}

Obviously $A_m(i,0)$ is $0$ unless $i=0$, in which case it is $1$.
If $n>0$ and $i=0$, then the number of walks is
$m A_m(1,n-1)$, since the first step must be to a vertex at distance $1$ from the starting vertex. If $i>0$ then we have the recurrence
$$
A_m(i,n)=A_m(i-1,n-1)+(m-1)A_m(i+1,n-1)
$$
for $i \geq 1$ and $n \geq 1$, since any vertex whose distance is $i$ from the starting point has exactly one
neighbor whose distance is $i-1$ and exactly $m-1$ neighbors whose distance is $i+1$.

With the same effort, we can solve a slightly more general problem.

\section*{A more general recurrence and weighted Dyck paths}

Let's consider, for any constants $c_1,c_2,c_3$, the solution of
the recurrence
\begin{align*}
A(i,n)&=c_1 A(i-1,n-1)+c_2 A(i+1,n-1), \quad i \geq 1, n \geq 1; \\
A(0,n)&= c_3 A(1,n-1), \quad n \geq 1; \\
A(i,0)&=\delta_{i,0}.
\end{align*}
This is a slight variation on the generic recurrence for Dyck paths. 
These are walks in the two-dimensional lattice with steps
\begin{align*}
	U&: (x,y) \to (x+1,y+1) \\
	D&: (x,y) \to (x+1,y-1)
\end{align*}
that start at the origin, end at the point $(n,i)$, and never venture below the $x$-axis.

Let's define the \emph{weight} of such a walk to be
$c_1^{\# U} c_2^{\# D} t^{\#U+\#D}$ and the \emph{poids} to be
$$
c_1^{\# U} c_2^\text{\#$D$'s not touching the $x$-axis} c_3^\text{\#$D$'s touching the $x$-axis} t^{\#U+\#D},
$$
or more succinctly,  
$$
c_1^{\# U} c_2^{\# D} (c_3/c_2)^\text{\# irreducible components} t^{\#U+\#D},
$$
where an irreducible component is a portion of a walk that starts and ends on the
$x$-axis, but otherwise is strictly above it.

It is easy to see that $A(i,n)t^n$ is the weight enumerator, according to \emph{poids}, of
the set of paths that terminate at $(n,i)$,
and $d_i(t):=\sum_{n=0}^\infty A(i,n)t^n$ is the
\emph{poids} enumerator of \emph{all} walks that end at a point with $y=i$.

Let's first compute $d_0(t)$, the \emph{poids} enumerator of 
walks that end on the $x$-axis (i.e for which $\#U=\#D$).
This is a minor variation on the usual way of counting Dyck paths.
Let $a(t)$ be the weight enumerator for \emph{all} Dyck paths according
to \emph{weight}, and let $b(t)$ be the weight enumerator (also according
to \emph{weight}) of \emph{irreducible} paths. Obviously, we have
$$
a(t)=\frac{1}{1-b(t)}, \quad b(t)=c_1c_2t^2 a(t).
$$
That leads to the quadratic equation
$$
(c_1c_2t^2)a(t)^2 -a(t)+ 1 = 0
$$
whose solution is
$$
a(t)= \frac{ 1 - \sqrt{1-4c_1c_2t^2} }{2c_1c_2t^2},
$$
and hence $b(t)=c_1c_2t^2a(t)$ equals
$$
b(t)= \frac{ 1 - \sqrt{1-4c_1c_2t^2} }{2}.
$$
Now let $c(t)$ be the weight enumerator of irreducible paths
according to \emph{poids}. Of course
$$
c(t)=\frac{c_3}{c_2} b(t), 
$$
since the last step now contributes a factor of $c_3$ rather than
$c_2$, from which we get
$$
c(t)= \frac{ c_3(1 - \sqrt{1-4c_1c_2t^2}) }{2c_2},
$$
and consequently, $d(t)=d_0(t)$, the weight enumerator of all Dyck paths according
to the weight \emph{poids} is
$$
d(t)= \frac{1}{1-c(t)}=
\frac{2c_2}{2c_2-c_3+c_3\sqrt{1-4c_1c_2t^2}}.
$$
By plugging in $c_1=1, c_2=m-1, c_3=m$, we get Paul Boddington's expression.

Finally, to get the weight enumerator $d_i(t)$ according to \emph{poids} of \emph{all} walks that
wind up at a point with $y=i$, note that any such walk
can be uniquely decomposed into $W_0 U W_1 U W_2 U W_3 \cdots U W_i$ where there are $i$ $U$'s
and $W_0, \dots, W_i$ are ordinary Dyck paths. Their \emph{poids} is
the \emph{poids} of $W_0$ times the product of the weights of $W_1, \dots , W_i$, times 
$c_1^i t^i$ (to account for the $i$ $U$'s). 
Hence
$$
d_i(t)= d(t) (c_1 t a(t))^i=
\frac{2c_2}{2c_2-c_3+c_3\sqrt{1-4c_1c_2t^2}} \cdot
\left( \frac{ 1 - \sqrt{1-4c_1c_2t^2} }{2c_2t} \right )^i.
$$
Plugging in $c_1=1,c_2=m-1,c_3=m$ yields
$$
f_m^{(i)}(t)=
\frac{2(m-1)}{m-2+m\sqrt{1-4(m-1)t^2}} \cdot
\left( \frac{ 1 - \sqrt{1-4(m-1)t^2} }{2(m-1)t} \right )^i.
$$

\end{document}